\documentclass[3p,times]{article}

\usepackage{graphicx, caption, subcaption}
\usepackage{amssymb}
\usepackage{amsmath}
\usepackage{amsthm}
\usepackage[final]{showkeys}
\usepackage{microtype}
\usepackage{color}
\usepackage{graphicx}
\usepackage{ulem}
\usepackage[figuresright]{rotating}
\usepackage[hmargin=3cm,vmargin={3.5cm,4cm}]{geometry}
\usepackage{bm}

\newtheorem{te}{Theorem}[section]
\newtheorem{example}{Example}

\newtheorem{os}[te]{Remark}
\newtheorem{prop}[te]{Proposition}

\numberwithin{equation}{section}
\title{A note on exact results for Burgers-like equations involving Laguerre derivatives.}

\author{Giuseppe Dattoli$^1$, Riccardo Droghei$^2$, Roberto Garra$^3$\\	
	$^1$ ENEA - Frascati Research Center, Via Enrico Fermi 45, 00044, \\Frascati, Rome, Italy  \\
	$^2$ CNR, ISMAR, Italy\\
 $^3$Università Telematica Internazionale, Section of Mathematics,  Italy\\
}

\date{\today}
\begin{document}
	
	\maketitle

	\abstract{In this note, we consider some Burgers-like equations involving Laguerre derivatives and demonstrate that it is possible to construct specific exact solutions using separation of variables. We prove that a general scheme exists for constructing exact solutions for these Burgers-like equations, extending to more general cases, including nonlinear time-fractional equations. Exact solutions can also be obtained for KdV-like equations involving Laguerre derivatives.
	We finally consider a particular class of Burgers equations with variable coefficients whose solution can be obtained similarly. }\\
	

\section{Introduction}

The Laguerre derivative has been introduced by A. Torre and one of the present authors (G. D.) in \cite{11}, within the context of a program aimed at establishing the theoretical framework of \textit{Monomiality} \cite{mon}.
This restatement of the theory of classical polynomials has a long history, tracing back to the calculus of differences, which has evolved significantly since the XVIII-th century, eventually merging with more recent umbral techniques (we refer for example to the recent monography \cite{umbr} and the references therein). Then, several research papers have been devoted to analysing Laguerre differential equations and related applications and special functions, we refer to the recent survey \cite{ricci} by P.E. Ricci.

The Laguerre polynomials are closely related to the properties of Bessel functions, as evidenced by the relevant generating function writes \cite{Andrews}
\begin{equation}
\sum_{n=0}^{\infty}\frac{t^n}{n!} L_n(x,y)= e^{yt}C_0(xt),
\end{equation}
where
\begin{equation}
L_n(x,y) = n!\sum_{r=0}^n\frac{(-1)^r x^r y^{n-r}}{r!^2(n-r)!}
\end{equation}
and
\begin{equation}
    C_0(x) = J_0(2\sqrt{x}) = \sum_{r=0}^\infty \frac{(-1)^r x^r}{r!^2},
\end{equation}
is the 0-th order Tricomi function or Laguerre exponential function \cite{ricci}. Even though it is a particular case of the ordinary cylindrical Bessel of the first kind, its properties are sufficiently interesting not to be considered a simple corollary of this family of functions. They have been indeed recognized as a generalization of the exponential functions \cite{dattoliricci}, satisfying the eigenvalue equation
\begin{equation}\label{eig}
	_L\widehat{D}_x C_0(\lambda x) = \lambda C_0(\lambda x),
\end{equation}
where
$$_L\widehat{D}_x := -\frac{d}{dx} x\frac{d}{dx} .$$
The operator $_L\widehat{D}_x$  is known as the Laguerre derivative, albeit, as underscored below, it is a particular case of the Bessel operator. This definition is however justified by the identity
\begin{equation}
    _L\widehat{D}_x L_n(x,y)= nL_{n-1}(x,y),
\end{equation}
which has been the starting point of recent research on the extension of the non-linear Burgers equation \cite{sachdev} to Laguerre type diffusive equations \cite{silvia}.\\
The motivation for this type of research has been manifold:\\
a)	The search for new non-linear equations associated with Laguerre/Bessel (LB) operators,\\
b)	The extension of the Hopf-Cole transformation to polynomials, not recognized as belonging to the Appèl family.\\
These polynomials are characterized by the generating function
\begin{equation}
    \sum_{n=0}^{\infty}\frac{t^n}{n!}a_n(x) = A(t) e^{xt},
\end{equation}
which, because the exponential function is an eigenfunction of the operator derivative, eventually yields 
\begin{equation}
    \sum_{n=0}^{\infty}\frac{t^n}{n!}a_n(x) = A(\partial_x) e^{xt},
\end{equation}
Finally expanding the exponential, we end up with the following operational definition of the Appèl polynomials \cite{babusci}
\begin{equation} \label{16a}
    a_n(x) = A(\partial_x) x^n.
\end{equation}
If we extend the same reasoning to eq. (1), we can conclude that
\begin{equation}\label{16b}
    L_n(x,y) = e^{y _L\widehat{D}_x}l_n(x),
\end{equation}
with
$$l_n(x) = \frac{(-x)^n}{n!}.$$

By comparing eqs.\eqref{16a} and \eqref{16b}, we can envisage the following correspondence
\begin{align}
 \partial_x \leftrightarrow  _L\widehat{D}_x,\\
    x^n \leftrightarrow l_n(x)
\end{align}
and 
\begin{align}\label{nonf}
\partial_x^m x^n = \frac{n!}{(n-m)!}x^{n-m}, \quad n\geq m, \\
\label{sec} _L\widehat{D}_x^m l_n(x) = \frac{n!}{(n-m)!}l_{n-m}(x) \quad n\geq m.
\end{align}
The identity in \eqref{sec} can be proven by noting that the LB derivative operator satisfies the property
\begin{equation}
    _L\widehat{D}_x^m  = (-1)^m \partial_x^m x^m \partial_x^m.
\end{equation}

\smallskip

In the recent papers \cite{silvia} and \cite{tomo}, the authors studied nonlinear equations involving Laguerre derivatives. Taking inspiration from the recent results discussed by Wazwaz in \cite{waz}, here we consider different classes of Burgers-like equations involving the Laguerre derivatives instead of the first-order time-derivative, showing that it is simple to obtain a solution by using the 0-th order Tricomi function. Then, we show that this simple construction of exact solutions for Burgers-like equations can be used for more general integro-differential Burgers-like equations, including the time-fractional one. We underline that the solution of nonlinear time-fractional equations is generally non-trivial. A similar scheme can be simply adapted to solve also Korteweg-de Vries-like equations, obtaining a different dispersion relation. \\
We finally consider another class of particular Burgers-like equations involving Laguerre derivatives and variable coefficients. In the conclusions, we suggest some open problems that can be considered for future research about the generalization of Laguerre derivatives and polynomials.

\section{Burgers-like equations involving Laguerre derivatives}
 
In the paper \cite{waz}, the author considered different Burgers-like equations of the form 
\begin{equation}
	u_t+Vu_x = \delta u_{xx},
\end{equation}
i.e. an advection-dissipation-like equation, where $V = V(u,u_x, u_{xx},\dots)$ is an arbitrary function. In \cite{waz} some interesting particular cases were considered, admitting simple travelling wave
type solutions. \\
The simplest case is the following one
\begin{equation}\label{eq1}
	u_t +\bigg(\frac{2u_x}{u}\bigg)u_x-u_{xx}=0
\end{equation}
whose solution is given by 
\begin{equation}
	u(x,t) = Re^{kx-rt},
\end{equation}
under the condition that $r = k^2$, while $R$ is an arbitrary real number. In \cite{waz} also other cases have been considered of Burger-like equations of this kind with a less straightforward form of
$V$. \\
According to the analogy established in the previous section, we use an extension of the previous equation in which the ordinary time derivative is replaced by the LB operator.\\
Let us consider the following generalized Burgers-like equation
\begin{equation}\label{eq2}
		-_L\widehat{D}_t u +\bigg(\frac{2u_x}{u}\bigg)u_x-u_{xx}=0.
\end{equation}
Inspired by the travelling wave solution, we can construct a simple solution using separation of variables. The equation \eqref{eq2}
admits a solution of the form
\begin{equation}\label{c0}
	u(x,t) = R e^{kx}C_0(rt),
\end{equation}
under the condition that $r = k^2$.
The travelling wave-like solution \eqref{c0} has the same structure of the generating function (1.1), thus suggesting that the solution of this type of non-linear equations can be also written in terms of Laguerre polynomials.
\\

This construction can be generalized, leading to the following result.
\begin{prop}
	The Burgers-type equation 
	\begin{equation}\label{eq3}
		\widehat{O}_t u +\bigg(\frac{2u_x}{u}\bigg)u_x-u_{xx}=0,
	\end{equation}
involving a general linear integro-differential operator $\widehat{O}_t$, admits a solution of the form 
\begin{equation}
u(x,t) = R e^{kx}f(-rt),
\end{equation}
where $f(\cdot)$ is such that 
$$\widehat{O}_t f(-rt) = -r f(-rt)$$
and $r = k^2$
\end{prop}
\begin{example}
	For example, we can consider the time-fractional Burgers-like equation
	\begin{equation}\label{ex}
		\widehat{D}^\alpha_t u +\bigg(\frac{2u_x}{u}\bigg)u_x-u_{xx}=0,
		\quad \alpha \in (0,1),
	\end{equation}
where we denoted by $\widehat{D}^\alpha_t$ the so-called Caputo fractional derivative (see \cite{kilbas}).
Recalling that the one-parameter Mittag-Leffler function 
\begin{equation}
	E_\alpha(-rt^\alpha) = \sum_{k=0}^\infty \frac{(-rt^\alpha)^k}{\Gamma(\alpha k +1)},
\end{equation}
is an eigenfunction of the Caputo derivative of order $\alpha$, we have that an exact solution for the equation \eqref{ex} is given by 
\begin{equation}
	u(x,t)  = R e^{kx}	E_\alpha(-rt^\alpha).	
\end{equation}
Observe that for $\alpha = 1$, we recover the travelling wave solution obtained in \cite{waz}. We should underline that in a few cases, it is possible to obtain explicit analytic solutions for nonlinear fractional Burgers-type equations, we refer for example to \cite{bur}. 
\end{example}
\begin{example}\label{HBexample}
A fractional generalization of the Laguerre derivative, analyzed in the previous section, can be obtained by introducing the 3-parameter hyper-Bessel type operator (see \cite{droghei2021},\cite{droghei2023})

\begin{equation}\label{hb}
_t\hat{D}_{\alpha,\beta,\nu}f(x)=x^{\alpha-\nu}\frac{d^\beta}{dx^\beta}\left(x^\nu\frac{d^\alpha}{dx^\alpha}f(x)\right)=f(x),
\end{equation}
involving two fractional derivatives in the sense of Caputo of orders $\alpha,\beta\in\left(0,1\right)$. Where $$f(x)=\mathcal{W}_{\alpha,\beta,\nu}(x^\beta)=\sum_{k=0}^\infty\prod_{i=1}^k\frac{\Gamma(\beta i+1-\alpha)}{\Gamma(\beta i+1)}\frac{x^{\beta k}}{\Gamma(\beta k+1-\alpha+\nu)}.$$
Let us consider the following generalized Burgers-like equation
\begin{equation}\label{eq2.2}
		_t\hat{D}_{\alpha,\beta,\nu} u +\bigg(\frac{2u_x}{u}\bigg)u_x-u_{xx}=0.
\end{equation}
As in eq.(\ref{eq2}), we can construct a simple solution by the separation of variables. The equation \eqref{eq2.2}
admits a solution of the form
\begin{equation}
	u(x,t) = R e^{kx}\mathcal{W}_{\alpha,\beta,\nu}(-r t^\beta),
\end{equation}
under the condition that $r = k^2$.\\
Obviously, for $\alpha=\beta=\nu=1$ the equation (\ref{eq2.2}) is reduced to equation (\ref{eq2}) and $$\mathcal{W}_{1,1,1}(t)=C_0(-t).$$
\end{example}
A similar construction can be simply adapted to the other cases considered in \cite{waz}.\\
For example, we have the nonlinear Burgers-like equation
\begin{equation}\label{eq5}
	\widehat{O}_t u +\bigg(u-\left(\frac{u_x^2+u_{xxx}}{u_{xx}}\right)\bigg)u_x+2u_{xx}=0,
\end{equation}
involving a general linear integro-differential operator $\widehat{O}_t$, admits a solution of the form 
\begin{equation}
	u(x,t) = R e^{kx}f(-rt),
\end{equation}
where $f(\cdot)$ is such that 
$$\widehat{O}_t f(-rt) = -r f(-rt)$$
and $r = k^2$.
We have a general result that can be adopted to solve a wide class of nonlinear integro-differential equations. This means that, for example, we have a scheme to construct exact solutions for Burgers-like equations involving higher-order Laguerre or Bessel derivatives. 
\begin{os}
We observe that this construction of simple exact solutions can be adapted also to the nonlinear equations 
\begin{equation}\label{eqn}
		-_L\widehat{D}_t u +\bigg(\frac{2u_x}{u}\bigg)^n u_x-u_{xx}=0, \quad n \geq 2.
\end{equation}
This general class of equations admits a solution of the form \eqref{c0} under the algebraic condition 
\begin{equation}
    r  = 2 k^{n+1}-k^2.
    \end{equation}
\end{os}

\bigskip

We now underline that this construction can be simply adapted to obtain exact solutions for Korteweg-de Vries-like equations \cite{waz2}.\\
Let us consider the following KdV-like equation 
\begin{equation}\label{eqkd}
		-_L\widehat{D}_t u +\bigg(\frac{2u_{xx}}{u}\bigg)u_x-u_{xxx}=0.
\end{equation}
It is simple to prove that the equation \eqref{eqkd}
admits a solution of the form
\begin{equation}
	u(x,t) = R e^{kx}C_0(rt),
\end{equation}
under the different dispersion relation, that is $r = k^3$.\\
Therefore, in analogy with the previous analysis, we can easily construct exact solutions for KdV-like equations
\begin{equation}
		\widehat{O}_t u +\bigg(\frac{2u_{xx}}{u}\bigg)u_x-u_{xxx}=0,
\end{equation}
involving a general linear integro-differential operator $\widehat{O}_t$, including the time-fractional KdV-like equations.\\
A similar scheme can be then adapted to solve other interesting nonlinear equations.

\section{Burgers equations with variable coefficients involving Laguerre derivatives}

We have shown that it is possible to construct particular solutions for nonlinear Burgers-like equations involving Laguerre derivatives. \\
Here we prove that a similar scheme can be simply adopted to solve Burgers equation with variable coefficients of the form 
\begin{equation}\label{beq}
		-_L\widehat{D}_t u +a(x)uu_x-b(t)u_{xx}+r u=0.
\end{equation}
If we consider a simple (but reasonable for the applications) choice of the time-variable coefficients $a(t)$ and $b(t)$.\\
Indeed, by taking $a(x) = k e^{-kx}$ and $b(t) = C_0(rt)$, we have that the Burgers-like solution \eqref{beq} admits a solution of the form 
\begin{equation}
	u(x,t) = e^{kx}C_0(rt).
\end{equation}

Also in this case we have a more general formulation
\begin{prop}
	The Burgers-type equation 
	\begin{equation}\label{eq3bis}
		\widehat{O}_t u +k e^{-kx}uu_x-b(t)u_{xx}+r u=0,
	\end{equation}
	involving a general linear integro-differential operator $\widehat{O}_t$, admits a solution of the form 
	\begin{equation}
		u(x,t) =  e^{kx}f(-rt),
	\end{equation}
	where $f(\cdot)$ is such that 
	$$\widehat{O}_t f(-rt) = -r f(-rt)$$
	and $b(t) = f(-rt)$.
\end{prop}

\section{Final comments and conclusions}

In this paper, we consider a general class of Burgers-like equations involving Laguerre derivatives, which admit exact solutions through the Laguerre exponential function. This scheme to construct exact solutions can be simply adapted to more general cases. Moreover, we have shown that this method can be applied to solve also KdV-like equations.\\

We emphasize that our focus was on nonlinear equations involving Laguerre time derivatives, although it would be interesting to explore Burgers-like or nonlinear diffusive-like equations involving Laguerre space derivatives.

Another interesting open problem that will be the object of further analysis is the study of the linear and non-linear equation involving fractional Laguerre-type derivatives. 
This fractional operator can be treated by starting from the equation \eqref{nonf} which can be extended formally to the fractional case too. We propose therefore the following generalization
\begin{equation}\label{4.1}
    _L\widehat{D}_x^\alpha l_n(x)  = \frac{n!}{\Gamma(n-\alpha+1)}l_{n-\alpha}(x)
\end{equation}
where $\alpha$ is real.
The last equation derives from a formal generalization of eq. \eqref{sec} to the non-integer case and provides a first hint to the generalization of LB derivatives to the fractional case.\\
Considering, furthermore, the Example \ref{HBexample} we can define the following fractional extension of the Laguerre-type derivative $_L\widehat{D}_x^{\alpha,\beta,\nu}=(-1)^\beta x^{\alpha-\nu}\frac{d^\beta}{dx^\beta}\left(x^\nu\frac{d^\alpha}{dx^\alpha}\right)$.
Indeed, we get 
\begin{equation}
    _L\widehat{D}_x^{\alpha,\beta,\nu} l^{(\alpha,\nu)}_n(x)  = \frac{n!}{\Gamma(n-\alpha+1)}l^{(\alpha,\nu)}_{n-\beta}(x)
\end{equation}
that becomes (\ref{4.1}) for $\beta=\nu=\alpha$, with 
\begin{equation}
    l^{(\alpha,\nu)}_n(x)=\frac{(-x)^n}{\Gamma(n-\alpha+\nu+1)},
\end{equation}
or
\begin{equation}
    _L\widehat{D}_x^{\alpha,\beta,\nu} \Tilde l^{(\alpha,\beta,\nu)}_n(x)  =\Tilde l^{(\alpha,\beta,\nu)}_{n-1}(x)
\end{equation}
with
\begin{equation}
    \Tilde l^{(\alpha,\beta,\nu)}_n(x)=\left(\prod_{i=1}^n\frac{\Gamma(\beta i+1-\alpha)}{\Gamma(\beta i+1)}\right) \frac{(-x)^{\beta n}}{\Gamma(\beta n-\alpha+\nu+1)}
\end{equation}
where $\alpha, \,\beta\in (0,1)$ and $\nu>0$.

\end{document}